\documentclass[reqno,12pt]{amsart}

\usepackage{amsmath}
\usepackage{amssymb}
\usepackage{amsthm}
\usepackage[english]{babel}
\newtheorem{theorem}{Theorem}

\newtheorem{proposition}{Proposition}

\newtheorem{corollary}{Corollary}

\begin{document}

\title[Estimates for the concentration functions]
{Estimates for the concentration functions in the Littlewood--Offord problem}

\author[Yu.S. Eliseeva]{Yulia S. Eliseeva}
\author[F. ~G\"otze]{Friedrich G\"otze}
\author[A.Yu. Zaitsev]{Andrei Yu. Zaitsev}

\email{pochta106@yandex.ru}
\address{St.~Petersburg State University\bigskip}
\email{goetze@math.uni-bielefeld.de}
\address{Fakult\"at f\"ur Mathematik,\newline\indent
Universit\"at Bielefeld, Postfach 100131,\newline\indent D-33501 Bielefeld,
Germany\bigskip}
\email{zaitsev@pdmi.ras.ru}
\address{St.~Petersburg Department of Steklov Mathematical Institute
\newline\indent
Fontanka 27, St.~Petersburg 191023, Russia\newline\indent
and St.~Petersburg State University}

\begin{abstract}{Let $X,X_1,\ldots,X_n$ be independent identically distributed random variables.
In this paper we study the behavior of the concentration functions
of the weighted sums $\sum_{k=1}^{n}a_k X_k$ with respect to the
arithmetic structure of coefficients~$a_k$. Such concentration
results recently became important in connection with
investigations about singular values of random matrices. In this
paper we formulate and prove some refinements of a result of
Vershynin (2014).}
\end{abstract}

\keywords {concentration functions, inequalities,
the Littlewood--Offord problem, sums of independent random variables}

\subjclass {Primary 60F05; secondary 60E15, 60G50}

\maketitle

\section{Introduction}

This paper is an extended and modified version of preprint
\cite{23}.

Let $X,X_1,\ldots,X_n$ be independent identically distributed (i.i.d.) random variables
with common distribution $F=\mathcal L(X)$.
The L\'evy concentration function of a random variable $X$ is defined by the equality
$$Q(F,\lambda)=\sup_{x\in\mathbf{R}}F\{[x,x+\lambda]\}, \quad \lambda>0.$$
 Let $a=(a_1,\ldots,a_n)\in \mathbf{R}^n$, $a\ne0$.
In this paper we study the behavior of the concentration functions
of the weighted sums $S_a=\sum_{k=1}^{n}a_k X_k$ with respect to
the arithmetic structure of coefficients~$a_k$. Refined
concentration results for these weighted sums play an important
role in the study of singular values of random matrices (see, for
instance, Nguyen and Vu \cite{Nguyen and Vu}, Rudelson and
Vershynin \cite{Rudelson and Vershynin08, Rudelson and Vershynin},
Tao and Vu \cite{Tao and Vu, Tao and Vu2}, Vershynin
\cite{Vershynin}). In this context the problem is referred to as
the Littlewood--Offord problem (see also \cite{29, 27, 28}).

In the sequel, let $F_a$ denote the distribution of the sum
$S_a$, and let $G$ be the distribution of the symmetrized random variable
$\widetilde{X}=X_1-X_2$.
Let \begin{equation} \label{0}M(\tau)=\tau^{-2}\int_{|x|\leq\tau}x^2
\,G\{dx\}+\int_{|x|>\tau}G\{dx\}=\mathbf{E}
\min\big\{{\widetilde{X}^2}/{\tau^2},1\big\}, \quad \tau>0.\end{equation}
 The symbol $c$ will be used for absolute positive constants. Note that $c$ can be different in different (or even in the same) formulas.
We will write $A\ll B$ if $A\leq c B$. Also we will write $A\asymp
B$ if $A\ll B$ and $B\ll A$. For~${x=(x_1,\dots,x_n )\in\mathbf R^n}$
we will denote $\|x\|^2= x_1^2+\dots +x_n^2$
and
$\|x\|_\infty= \max_j|x_j|$.

The elementary properties of concentration functions are well studied (see, for instance,
 \cite{Arak and Zaitsev, Hengartner and Theodorescu,
Petrov}). In particular, it is obvious that $Q(F,\mu)\le
(1+\lfloor \mu/\lambda\rfloor)\,Q(F,\lambda)$, for any
$\mu,\lambda>0$, where $\lfloor x\rfloor$ is the integer part of a
number~$x$. Hence,
\begin{equation}\label{8a} Q(F,c\lambda)\asymp\,Q(F,\lambda)
\end{equation}
and\begin{equation}\label{8j}
 \hbox{if }Q(F,\lambda)\ll B,\hbox{ then }Q(F,\mu)\ll B\,(1+\mu/\lambda).
\end{equation}

The problem of estimating the concentration function of weighted
sums $S_a$ under different conditions on the vector $a\in
\mathbf{R}^n$ and distributions of summands has been studied in
\cite{Friedland and Sodin, Nguyen and Vu,  Rudelson and
Vershynin08, Rudelson and Vershynin, Tao and Vu, Tao and Vu2,
Vershynin}. Eliseeva and Zaitsev \cite{Eliseeva and Zaitsev}, see
also \cite{25}, have obtained some improvements of the results
\cite{Friedland and Sodin} and~\cite{Rudelson and Vershynin}.  In
this paper we formulate and prove similar refinements of a result
of Vershynin~\cite{Vershynin}.

Note that a connection of the rate of decay of the concentration
functions of sums  with the arithmetic structure of supports of
distributions of independent random variables was found by Arak
\cite{Arak}, see also \cite{Arak and Zaitsev}, long before the
appearance of the papers \cite{Friedland and Sodin, Nguyen and Vu,
Rudelson and Vershynin08, Rudelson and Vershynin, Tao and Vu, Tao
and Vu2, Vershynin} in which a similar connection was considered
in a particular case of distributions from the Littlewood--Offord
problem. The authors of the present paper are going to discuss
these connections in a separate publication.
\bigskip

Denote  $ \log_+(x)$ = $\max\{ 0,\log x\}$. The result of
Vershynin \cite{Vershynin}, related to the Littlewood--Offord
problem, is formulated as follows.

\begin{proposition}\label{thV} Let $X, X_1,\ldots,X_n$ be i.i.d. random variables
and
$a=(a_1,\ldots,a_n)\in \mathbf{R}^n$
 with $\|a\|=1$.  Assume that there exist positive numbers
 $\tau, p, K, L,D$ such that $Q(\mathcal{L}(X),\tau)\leq 1-p$,
 $\mathbf{E} \,\left|X\right|\le K$,  and
 \begin{equation} \label{5bk}
\|\,t a-m\|\geq L\sqrt{\log_+(t/L)}\ \hbox{ for all $m\in \mathbf Z^n$  and \
$t\in(0,D]$}.
\end{equation}
  If $L^2\geq 1/{p}$,  then
\begin{equation} \label{5bff} Q\Big(F_a,\cfrac{1}{D}\Big)\le \frac{C\,L}{D},\end{equation}
where the quantity $C$ depends on $\tau, p, K$ only.\end{proposition}

\begin{corollary}\label{c1a} Let the conditions of Proposition\/ $\ref{thV}$ be satisfied.
Then, for any $ \varepsilon\ge 0$,
\begin{equation} \label{7bb}Q(F_a,\varepsilon)
\ll C\,L\,\Big(\varepsilon+\cfrac{1}{D}\Big). \end{equation}
\end{corollary}
\medskip

 It is clear that if
 \begin{equation} \label{4bb}
0<D\le D(a)=D_{L}(a)=
\inf\Big\{t>0:\hbox{dist}(ta,\mathbf{Z}^n)<
 L\sqrt{\log_+(t/L)}\Big\},
\end{equation}where $$\hbox{dist}(ta,\mathbf{Z}^n)= \min_{m \in \mathbf{Z}^n}\|\,ta - m\|
=\Big(\sum_{k=1}^{n} \min_{m_k \in \mathbf{Z}} |\,ta_k -
m_k|^2\Big)^{1/2} ,
 $$then  condition \eqref{5bk} holds.
In Vershynin \cite{Vershynin}  the quantity~$D(a)$ is called
the least common denominator of the vector ${a\in\mathbf{R}^n}$ (see also
Rudelson and Vershynin \cite{Rudelson and Vershynin08} and \cite{Rudelson and
Vershynin} for similar definitions).

Note that for $|\,t|\leq 1/2\,\|a\|_{\infty}$ we have
\begin{equation}\label{4s}
\big(\hbox{dist}(ta,\mathbf{Z}^n)\big)^2 =
\sum_{k=1}^{n}|\,ta_k|^2= \|a\|^2t^2= t^2. \end{equation} Hence,
by definition, $D(a)> L$. Moreover, equality~\eqref{4s} implies
that $D(a)\ge {1}/{2\,\|a\|_{\infty}}$ (see
Vershynin~\cite{Vershynin}, Lemma~6.2).

Note that just the statement of Corollary~\ref{c1a} with $D=D(a)$
is formulated in~\cite{Vershynin}. Proposition \ref{thV} seems to
be more natural formulation which implies Corollary~\ref{c1a}
using relations \eqref{8j} and~\eqref{4bb}. Minimal $L$ satisfying
the conditions of Proposition \ref{thV} depends on $a$ and $D$,
and, generally, may be much larger than $p^{-1/2}$.

In the formulation of Proposition~\ref{thV},  w.l.o.g. we can replace assumption~\eqref{5bk} by the following:
\begin{equation} \label{5bh}
\|\,t a-m\|\geq f_L(t)\quad \hbox{ for all $m\in \mathbf Z^n$ and \
$t\in \Big[\cfrac{1}{2\,\|a\|_\infty},D\Big]$},
\end{equation}where
\begin{equation} \label{5cv}
f_L(t)= \begin{cases}\qquad t/6,& \hbox{ for }0<t< eL,\\ L\sqrt{\log(t/L)},&\hbox{ for }t\ge eL.
\end{cases}
\end{equation}
Note that equality~\eqref{4s} justifies why the assumption
$t\geq1/2\,\|a\|_{\infty}$ in condition~\eqref{5bh} is natural.
For $0<t< 1/2\,\|a\|_{\infty}$, inequality~\eqref{5bh} is
satisfied automatically.

Formally,  condition~\eqref{5bh} may be more restrictive than
condition~\eqref{5bk}.
 However, if condition~\eqref{5bk} is satisfied, but condition~\eqref{5bh} not,
 then inequality~\eqref{5bff} remains true by trivial reasons.

Indeed, if $t\ge eL$, then the validity of~\eqref{5bh} for such
a~$t$ follows from assumption~\eqref{5bk}. If $0<t< eL$ and there
exists an $m\in \mathbf Z^n$ such that $\|\,t a-m\|<t/6$, then,
denoting $k=\lfloor eL/t\rfloor+1$, we have $tk\ge eL$ and
$$\|\,tk a-km\|<tk/6\le 2eL/6<L\leq L\sqrt{\log_+(tk/L)}.$$ Since
$km\in\mathbf Z^n$,  we have $D\le D(a)\le tk\ll L$ and the
required inequality~\eqref{5bff} is a trivial consequence of
$Q(F_a, 1/D)\le1$.

Note that the exists a possibility that condition~\eqref{5bh} is
satisfied, but condition~\eqref{5bk} not, for some $t$ from the
interval $L<t<eL$. Then the bounds for concentration functions
from Proposition~\ref{thV} and Corollary~\ref{c1a} remain true.
This follows from Theorem~\ref{th1a} of the present paper.

The above arguments justify that the least common denominator
$D^*(a)$ should be defined as
 \begin{equation} \label{4bt}
D^*(a)=
\inf\Big\{t>0:\hbox{dist}(ta,\mathbf{Z}^n)<
f_L(t\|a\|) \Big\}.
\end{equation}
This definition will be also used below in the case when $\|a\|\ne1$.
Obviously,
\begin{equation} \label{4btr}
D^*(\lambda a)=D^*(a)/\lambda, \quad\hbox{ for any }\lambda>0,
\end{equation}
and equality~\eqref{4s} implies also that $D^*(a)\ge 1/2\,\|a\|_{\infty}$.

 \medskip

Now we formulate the main result of this paper.

\begin{theorem}\label{th1a} Let $X,X_1,\ldots,X_n$ be {i.i.d.} random variables. Let
$a=(a_1,\ldots,a_n)\in \mathbf{R}^n$
 with $\|a\|=1$. Assume that condition \eqref{5bh}
 is satisfied.
  If $L^2\geq 1/{M(1)}$, where the quantity $M(1)$ is defined by formula \eqref{0}, then
\begin{equation} \label{5bf} Q\Big(F_a,\cfrac{1}{D}\Big)\ll \frac{1}{D\sqrt{M(1)}}.\end{equation}
\end{theorem}
\medskip

Let us reformulate Theorem\/ $\ref{th1a}$ for arbitrary $a$, without assuming that $\|a\|=1$.

\begin{corollary}\label{c1b} {Let the conditions of Theorem\/ $\ref{th1a}$ be satisfied
without the assumption $\|a\|=1$ and with  condition~\eqref{5bh}
replaced by the condition
\begin{equation} \label{5bj}
\|\,t a-m\|\geq f_L(t\|a\|)\ \hbox{ for all $m\in \mathbf Z^n$ and
\ $t\in \Big[\cfrac{1}{2\,\|a\|_\infty},D\Big]$}.
\end{equation}
 If $L^2\geq 1/{M(1)}$, then
\begin{equation} \label{5bj3}
Q\Big(F_a,\cfrac{1}{D}\Big)\ll \cfrac{1}{\|a\|D\sqrt{M(1)}} \,.
\end{equation}}\end{corollary}
\medskip

The proofs of our  Theorem \ref{th1a} and Corollary \ref{c1b}
are similar to the proof of the main results of Eliseeva and Zaitsev \cite{Eliseeva and Zaitsev}.
They are in some sense more natural than the proofs in
 Vershynin~\cite{Vershynin}, since they do not use unnecessary assumptions like
$\mathbf{E} \,\left|X\right|\le K$. This is achieved by an
application of relation~\eqref{1b}. Our proof differs from the
arguments used in \cite{Friedland and Sodin, Rudelson and
Vershynin, Vershynin} since we rely on
 methods introduced by Ess\'een \cite{Esseen} (see the proof of Lemma~4 of Chapter~II  in
\cite{Petrov}).
\medskip

Now we reformulate Corollary \ref{c1b} for the random variables
${X_k}/{\tau}$, $\tau>0$.

\begin{corollary}\label{c2a} Let
$V_{a,\tau}=\mathcal{L}\big(\sum_{k=1}^{n}a_k {X_k}/{\tau}\big)$, $\tau>0$.
Then, under the conditions of Corollary~$\ref{c1b}$
with the condition~$L^2\geq 1/{M(1)}$ replaced by the condition
$L^2\geq 1/{M(\tau)}$, we have
\begin{equation}\label{4p}
Q\Big(V_{a,\tau},\cfrac{1}{D}\Big) = Q\Big(F_a,\cfrac{\tau}{D}\Big) \ll
\cfrac{1}{\|a\|D\sqrt{M(\tau)}}\,.
\end{equation}
In particular, if $\|a\|=1$, then
\begin{equation}\label{6p}Q\Big(F_a,\cfrac{\tau}{D}\Big) \ll
\cfrac{1}{D\sqrt{M(\tau)}}\,.\end{equation}
\end{corollary}
\medskip

 For the proof of Corollary \ref{c2a}, it suffices to use Corollary \ref{c1b} and relation \eqref{0}.
 \medskip

 It is evident that $M(\tau)\gg 1-Q(G,\tau)\geq 1-Q(F,\tau)\geq p$, under the conditions of
 Proposition \ref{thV}. Note that $M(\tau)$ may be essentially larger than~$p$.
 For example, $p$ may be equal to $0$, while
$M(\tau)>0$ for any non-degenerate distribution ${F=\mathcal
L(X)}$. Comparing the bounds \eqref{5bff} and \eqref{6p}, we see
that the factor~$L$ is replaced by the
factor~${1}/{\sqrt{M(\tau)}}\le L$ which can be essentially
smaller than~$L$ under the conditions of  Corollary~\ref{c2a}.
 Moreover, there is an unnecessary assumption $\mathbf{E} \,\left|X\right|\le K$  in the formulation of
Proposition \ref{thV}. Finally, the dependence of constants on the
distribution $\mathcal L(X)$ is stated explicitly, in inequalities
\eqref{5bf} and \eqref{5bj3}--\eqref{6p} the constants are
absolute,  in contrast with inequalities \eqref{5bff} and
\eqref{7bb}, where $C$ depends on $\tau, p$ and $K$ in a
non-explicit way. An improvement of Corollary \ref{c1a} is given
below in Theorem~\ref{th1b}.
\medskip

We recall now the well-known Kolmogorov--Rogozin inequality
\cite{Rogozin} (see \cite{Arak and Zaitsev, Hengartner and
Theodorescu, Petrov}).
\medskip

\begin{proposition}\label{thKR} Let $Y_1,\ldots,Y_n$ be independent random variables with the distributions $W_k=\mathcal{L}(Y_k)$.
Let $\lambda_1,\ldots,\lambda_n$
be positive numbers such that ${\lambda_k \leq \lambda}$, for $k=1,\ldots,n$. Then
\begin{equation} \label{2}
Q\Big(\mathcal{L}\Big(\sum_{k=1}^{n}Y_k\Big),\lambda\Big)\ll\lambda\,
\Big(\sum_{k=1}^{n}\lambda_k^2\,\big(1-Q(W_k,\lambda_k)\big)\Big)^{-1/2}.
\end{equation}
\end{proposition}

Ess\'een \cite{Esseen} (see \cite{Petrov}, Theorem 3 of Chapter III)
has improved this result. He has shown that the following statement is true.

\begin{proposition}\label{thE}Under the conditions of Proposition $\ref{thKR}$ we have
\begin{equation} \label{3}
Q\Big(\mathcal{L}\Big(\sum_{k=1}^{n}Y_k\Big),\lambda\Big)\ll\lambda
\,\Big(\sum_{k=1}^{n}\lambda_k^2\,
M_k(\lambda_k)\Big)^{-1/2},
\end{equation}
 where $M_k(\tau)=\mathbf{E}\,
\min\big\{{\widetilde{Y_k}^2}/{\tau^2},1\big\}$.\end{proposition}

Furthermore, improvements of \eqref{2} and \eqref{3} may be found
in \cite{Arak, Arak and Zaitsev, Bretagnolle, GZ1, GZ2, Kesten,
Miroshnikov and Rogozin} and~\cite{Nagaev and Hodzhabagyan}.
 \medskip

It is clear that Theorem \ref{th1a} is related to Proposition
\ref{thV} in a similar way as  Ess\'een's inequality~\eqref{3} is
related to the Kolmogorov--Rogozin inequality~\eqref{2}. In
addition, the dependence of $C$ on $\tau, p$ and $K$ in
 \eqref{5bff} and \eqref{7bb} is not written out explicitly.

If we consider a special case, where $D=1/2\,\|a\|_{\infty}$, then
no assumptions on the arithmetic structure of the vector $a$ are
made, and Corollary~\ref{c2a} implies the bound
\begin{equation}\label{4}
Q(F_a,\|a\|_{\infty}\,\tau)\ll
\cfrac{\|a\|_{\infty}}{\|a\|\sqrt{M(\tau)}}\,.
\end{equation}
This result follows from Ess\'een's inequality \eqref{3}
applied to the sum of non-identically distributed random variables $Y_k=a_kX_k$ with
$\lambda_k=a_k\,\tau$,
$\lambda=\|a\|_\infty\,\tau$.
For $a_1=a_2=\cdots=a_n=n^{-1/2}$, inequality \eqref{4} turns into
the well-known particular case of Proposition \ref{thE}:
\begin{equation}\label{4y}
Q(F^{*n},\tau)\ll \cfrac{1}{\sqrt{n\,M(\tau)}}\,.
\end{equation}
Inequality \eqref{4y} implies also the Kolmogorov--Rogozin
inequality for i.i.d. random variables:
$$Q(F^{*n},\tau)\ll \cfrac{1}{\sqrt{n\,(1-Q(F,\tau))}}\,.$$

Inequality \eqref{4} can not yield bound of better order than
$O(n^{-1/2})$, since the right-hand side of \eqref{4} is at least
$n^{-1/2}$. The results stated above are more interesting if $D$
is essentially larger than ~$1/2\,\|a\|_{\infty}$. In this case
one can expect the estimates of much smaller order than
$O(n^{-1/2})$. Such estimates of $Q(F_a,\lambda)$ are required to
study the distributions of eigenvalues of random matrices.

For
$0<D<1/2\,\|a\|_{\infty}$, the inequality
\begin{equation}\label{4m}
 Q\Big(F_a,\cfrac{\tau}{D}\Big) \ll
\cfrac{1}{\|a\|D\sqrt{M(\tau)}}
\end{equation}
holds under the conditions of Corollary \ref{c2a} too. In this
case it follows from \eqref{8j} and \eqref{4}.
 \medskip

Under the conditions of Corollary \ref{c2a}, there exist many possibilities to represent
a fixed $\varepsilon$ as $\varepsilon=\tau/D$ for an appication of inequality~\eqref{4p}.
Therefore, for a fixed $\varepsilon=\tau/D$ we can try to minimize the right-hand side of inequality~\eqref{4p} choosing an optimal~$D$.
This is possible, and the optimal bound is given
in the following Theorem~\ref{th1b}.

\begin{theorem}\label{th1b} Let
 the conditions of Corollary $\ref{c1b}$ be satisfied for $D\le
 D^*(a)$,
except the condition~$L^2\geq 1/{M(1)}$. Let $L^2> 1/{ P}$, where
$P={\mathbf P}(\widetilde{X}\neq 0)=\lim_{\tau\to0}M(\tau)$. Then
there exists a $\tau_0$ such that $L^2=1/{M(\tau_0)}$. Moreover,
the bound \begin{equation}\label{4pt2} Q\big(F_a,\varepsilon\big)
\ll \cfrac{1}{\|a\|D^*(a)\sqrt{M(\varepsilon\, D^*(a))}}
\end{equation}
is valid for $0<\varepsilon\le\varepsilon_0=\tau_0/D^*(a)$.
Furthermore, for $\varepsilon\ge\varepsilon_0$,
 the bound
\begin{equation}\label{4pr2}
Q\big(F_a,\varepsilon\big) \ll
\cfrac{\varepsilon L}{\varepsilon_0\,\|a\|D^*(a)}
\end{equation} holds.
\end{theorem}
\medskip

In the statement of Theorem~\ref{th1b}, the quantity~$\varepsilon$ can be arbitrarily small. If
$\varepsilon$ tends to zero, we obtain
\begin{equation}\label{4pp4}
Q(F_a,0)\ll \cfrac{1}{\|a\|D^*(a)\sqrt{P}}\,,
\end{equation}  if $L^2> 1/{P}$. Applying inequalities
\eqref{4pt2}--\eqref{4pp4},
one should take into account that, by \eqref{4btr},
$\|a\|D^*(a)=D^*(a/\|a\|)$.

Theorem\/ $\ref{th1b}$ follows easily from Corollary \ref{c2a}.
Indeed, denoting $\varepsilon=\tau/D$, we can rewrite inequality~\eqref{4p} as
\begin{equation}\label{4pp}
Q\big(F_a,\varepsilon\big) \ll
\cfrac{1}{\|a\|D\sqrt{M(\varepsilon\, D)}}\,.
\end{equation}
Inequality~\eqref{4pp} holds if $L^2\geq 1/{M(\varepsilon\, D)}$
and $0<D\le D^*(a)$. If $L^2\geq 1/{M(\varepsilon\, D^*(a))}$,
then the choice $D=D^*(a)$ is optimal in inequality~\eqref{4pp}
since $$D^2{M(\varepsilon\, D)}=\mathbf{E}
\min\big\{{\widetilde{X}^2}/{\varepsilon^2},D^2\big\}$$ is
increasing when $D$ increases. For the same reason, if $L^2<
1/{M(\varepsilon\, D^*(a))}$, then the optimal choice of $D$ in
inequality~\eqref{4pp} is given by the solution $D_0(\varepsilon)$
of the equation~${L^2= 1/{M(\varepsilon\, D)}}$. This solution
exists and is unique if $L^2> 1/{ P}$, since
 the function $M(\tau)$ is continuous and strictly decreasing if $M(\tau)<P$.
Moreover, it is clear that $M(\tau)\to0$ as $\tau\to\infty$. In
this case inequality~\eqref{4pp} turns into
\begin{equation}\label{4ppp}
Q\big(F_a,\varepsilon\big) \ll
\cfrac{L}{\|a\|D_0(\varepsilon)}.
\end{equation}
Moreover, choosing $\tau_0$ as the solution of the equation
$L^2=1/{M(\tau)}$, we see that inequality~\eqref{4pt2} is valid
for $0<\varepsilon\le\varepsilon_0=\tau_0/D^*(a)$. It is clear
that $D_0(\varepsilon_0)=D^*(a)$. Furthermore, for
$\varepsilon\ge\varepsilon_0$, we have
$$
M(\varepsilon\, D_0(\varepsilon))=M(\varepsilon_0\, D_0(\varepsilon_0))=L^{-2}
$$
and, hence, $\varepsilon\, D_0(\varepsilon)=\varepsilon_0\,
D_0(\varepsilon_0)$. Therefore, for $\varepsilon\ge\varepsilon_0$,
inequality~\eqref{4pr2} holds. The right-hand side of this
inequality with $\|a\|=1$ admits also representations
$$
\cfrac{\varepsilon L}{\varepsilon_0\,D^*(a)}=
\cfrac{L}{D_0(\varepsilon)}=
\cfrac{1}{D_0(\varepsilon)\sqrt{M(\varepsilon\,
D_0(\varepsilon))}}\,.
$$

Obviously, inequality \eqref{4pr2} could be derived from \eqref{4pp} with $\varepsilon=\varepsilon_0$
by an application of inequality~\eqref{8j}. On the other hand, for
$0<\varepsilon_1<\varepsilon\le\varepsilon_0$,
we could apply inequality~\eqref{8j} to inequality \eqref{4pt2} and obtain the bound
\begin{equation}\label{4ps}
Q\big(F_a,\varepsilon\big) \ll \frac\varepsilon{\varepsilon_1}\, Q\big(F_a,\varepsilon_1\big) \ll
\cfrac{\varepsilon}{\varepsilon_1\,\|a\|D^*(a)\sqrt{M(\varepsilon_1\, D^*(a))}}\,.
\end{equation}
However, inequality \eqref{4ps} is weaker than inequality \eqref{4pt2} since, evidently,
\begin{equation}\label{4pa}
\varepsilon^2M(\varepsilon\, \mu)=\mathbf{E}
\min\big\{{\widetilde{X}^2}/{\mu^2},\varepsilon^2\big\}\ge \mathbf{E}
\min\big\{{\widetilde{X}^2}/{\mu^2},\varepsilon_1^2\big\}=\varepsilon_1^2\,M(\varepsilon_1\,\mu),
\end{equation}
for any $\mu>0$.

Theorem~\ref{th1b} is an essential improvement of Corollary
\ref{c1a}. In particular, in contrast with inequality \eqref{7bb}
of Corollary~\ref{c1a}, for small $\varepsilon$, the right-hand
side of inequality~\eqref{4pt2} of Theorem~\ref{th1b} may be
decreasing as $\varepsilon$ decreases. Moreover, we have just
shown that the application of inequality~\eqref{8j} would lead to
a loss of precision. Recall that Corollary~\ref{c1a} could be
derived from Proposition~\ref{thV} with the help of
inequality~\eqref{8j}.

Consider a simple example. Let $X$
 be the random variable
taking values $0$ and $1$ with probabilities
\begin{equation}
{\bf P}\{ X=1\}=1-{\bf P}\{ X=0\}=p>0.\label{1p}
\end{equation}Then
\begin{equation}
{\bf P}\{ \widetilde X=\pm1\}=p(1-p),\quad{\bf P}\{
\widetilde X=0\}=1-2\,p(1-p),\label{11p}
\end{equation}
and the function $M(\tau)$ has the form
\begin{equation} \label{5hh}
M(\tau)= \begin{cases}\quad 2\,p(1-p),& \hbox{ for }
0<\tau< 1,\\ 2\,p(1-p)/\tau^2,&\hbox{ for }\tau\ge 1.
\end{cases}
\end{equation}
Assume for simplicity that $\|a\|=1$.
If $L^2>1/2\,p(1-p)$, then $\tau_0=L\sqrt{2\,p(1-p)}$ and,
for $\varepsilon\ge\varepsilon_0=L\sqrt{2\,p(1-p)}/D^*(a)$,
we have the bound
\begin{equation}
Q\big(F_a,\varepsilon\big) \ll \cfrac{\varepsilon
}{\sqrt{p(1-p)}}\,.\label{17p}
\end{equation}
The same bound~\eqref{17p} follows from inequality~\eqref{4pt2} of Theorem~\ref{th1b} for $1/D^*(a)\le\varepsilon\le\varepsilon_0$.
For $0<\varepsilon\le1/D^*(a)$, inequality~\eqref{4pt2} implies the bound
\begin{equation}
Q\big(F_a,\varepsilon\big) \ll \cfrac{1
}{D^*(a)\sqrt{p(1-p)}}\,.\label{18p}
\end{equation}
Thus, \begin{equation}
Q\big(F_a,\varepsilon\big) \ll\min\bigg\{
\cfrac{1 }{\sqrt{p(1-p)}}\Big(\varepsilon+\cfrac{1}{D^*(a)}\Big),\:1\bigg\},
\quad\hbox{for all }\varepsilon\ge0.\label{19p}
\end{equation}
Inequality~\eqref{19p} cannot be essentially improved. Consider,
for instance,
\begin{equation}\label{19py}
a=(s^{-1/2}, \ldots, s^{-1/2},0,\ldots,0)
 \end{equation} with
the first $s\le n$ coordinates equal to~$s^{-1/2}$ and
 the last $n-s$ coordinates equal to zero. In this case $D^*(a)\asymp s^{1/2}$,
the random variable $s^{1/2}S_a$ has binomial distribution with parameters $s$ and $p$,
and it is well-known that
\begin{equation}
Q\big(F_a,\varepsilon\big) \gg\min\bigg\{
\cfrac{1 }{\sqrt{p(1-p)}}\,\Big(\varepsilon+\cfrac{1}{\sqrt s}\Big),\:1\bigg\},
\quad\hbox{for all }\varepsilon\ge0.\label{13pp}
\end{equation}
Comparing the bounds \eqref{19p} and \eqref{13pp}, we see that
Theorem~\ref{th1b} provides the optimal order of
$Q\big(F_a,\varepsilon\big)$ for all possible values
of~$\varepsilon$. Moreover, the involved constant is absolute.

It may seem that the last example is reduced to a trivial case
$n=s$. This is not entirely true. Clearly, the value of
$Q\big(F_a,1\big)$ cannot be changed much for a small change of
the vector~$a$, defined in~\eqref{19py}, if the last $n-s$
coordinates of this vector are small in magnitude, but not zero.
The degree of smallness of the last $n-s$ coordinates can be
chosen so that inequalities \eqref{19p} and \eqref{13pp}
 remain true for $\varepsilon\gg s^{-1}$ and $D^*(a)\asymp s^{1/2}$.
 \medskip

For the sake of completeness, we give below a short proof of
inequality~\eqref{13pp}. It is easy to see that
$\hbox{Var}\big(S_a\big)=p(1-p)$. Hence, by Chebyshev's
inequality,
\begin{equation}
\mathbf{P}\big\{|S_a-\mathbf{E}\,S_a|<2\sqrt{p(1-p)}\big\}\ge3/4.\label{13r}
\end{equation}
The random variable $S_a$ takes values which are multiples
of~$s^{-1/2}$. Therefore, if $s\,p(1-p)\le1$, then
inequality~\eqref{13r} implies that $Q\big(F_a,0\big)\asymp1$ and
inequality~\eqref{13pp} is valid.

Assume now $s\,p(1-p)>1$.
If $0<\varepsilon\le4\sqrt{p(1-p)}$, then, using \eqref{8j}  and~\eqref{13r},
we obtain
\begin{equation}
3/4\le Q\big(F_a,4\sqrt{p(1-p)}\big)\ll
\varepsilon^{-1}\sqrt{p(1-p)}\,Q\big(F_a,\varepsilon\big),\label{15r}
\end{equation} and, hence,
\begin{equation}
Q\big(F_a,\varepsilon\big) \gg
\cfrac{\varepsilon}{\sqrt{p(1-p)}}\,.\label{2pp}
\end{equation}
It is clear that \eqref{8a}, \eqref{8j}  and~\eqref{2pp} imply that $Q\big(F_a,\varepsilon\big)\asymp1$,
for $\varepsilon\ge4\sqrt{p(1-p)}$.
Applying inequality~\eqref{2pp} for $\varepsilon=s^{-1/2}$
and using the lattice structure of the support of distribution~$F_a$,
we conclude that,
for $0\le\varepsilon<{s^{-1/2}}$,
\begin{equation}
Q\big(F_a,\varepsilon\big) \ge Q\big(F_a,0\big)\gg
\cfrac{1 }{\sqrt{s\,p(1-p)}}\,.\label{12pp}
\end{equation}
Thus, inequalities~\eqref{8a}, \eqref{8j}, \eqref{2pp} and~\eqref{12pp} imply \eqref{13pp}.
 \medskip

 The results of this paper are formulated for a fixed~$ L $.
It is clear that in their application one should try to choose an
optimal~$ L $, satisfying the assumptions and minimizing the
right-hand sides of inequalities, which give bounds for the
concentration functions. Recall that the least common denominator
$D^*(a)$ depends on~$ L $.

The quantity $\tau_0=\varepsilon_0\, D^*(a)$ (which is the
solution of the equation ${L^2=1/{M(\tau)}}$) may be interpreted
as a quantity depending on $L$ and on the
distribution~$\mathcal{L}(X)$. Moreover, comparing the bounds
\eqref{7bb} and \eqref{4pr2} for relatively large values
of~$\varepsilon$, we see that $\tau_0\to\infty$  as $L\to\infty$.
Therefore, the factor $L/\tau_0$ is much smaller than~$L$ for
large values of~$L$. In particular, in the above example we have
$\tau_0=L\sqrt{2\,p(1-p)}$.

Another example would be a symmetric stable distribution with
parameter~$\alpha$, $0<\alpha<2$. In this case the characteristic
function $\widehat{F}(t)= \mathbf{E}\,\exp(itX)$ has the
form~$\widehat{F}(t)= \exp(-c\,|t|^\alpha)$. It could be shown
that then $\tau_0$ behaves as $L^{2/\alpha}$ as $L\to\infty$.

Inequality \eqref{17p} can be rewritten in the form
\begin{equation}
Q\big(F_a,\varepsilon\big) \ll
\cfrac{\varepsilon }{\sigma},\quad\hbox{for }\varepsilon\ge\varepsilon_0,\label{117p}
\end{equation}
where $\sigma^2=\hbox{Var}(X)$. It is clear that a similar
situation occurs for any random variable~$X$ with finite variance.

In particular, inequality \eqref{117p} is obviously satisfied for
all $\varepsilon\ge0$, if $\|a\|=1$ and $X$ has a Gaussian
distribution with $\hbox{Var\;}(X)=\sigma^2$. The order of this
inequality is optimal for $0\le\varepsilon\le\sigma$. In this
specific case, the relation
$$
\frac1{\sqrt{M(\tau)}}\asymp1+\frac{\tau}{\sigma},
$$
holds, for any $\tau>0$. Together with Theorem~\ref{th1b} for
$\|a\|=1$, it implies easily that
\begin{equation}\label{169p} Q\big(F_a,\varepsilon\big)
\ll\cfrac{\varepsilon }{\sigma}\quad \hbox{for }\varepsilon\ge
\cfrac{\sigma}{D^*(a)}.
\end{equation}
This gives the correct dependence of the concentration function on
 $\sigma$ if $\sigma/D^*(a)\le\varepsilon\le\sigma$. The same order
 of the bound can not be achieved with the help of
inequality~\eqref{7bb}.  It is impossible to derive
inequality~\eqref{169p} for small $\varepsilon$ from
Theorem~\ref{th1b}. This is due to the fact that the distribution
${F=\mathcal L(X)}$ is arbitrary in Theorem~\ref{th1b}, and the
concentration function $Q\big(F_a,\varepsilon\big)$ can not tend
to zero as $\varepsilon\to0$ (see~\eqref{13pp}).
\bigskip

\section{Proofs}

We will use the classical Ess\'een inequalities (\cite{Esseen66}, see also
\cite{Hengartner and Theodorescu} and \cite{Petrov}):
\begin{equation} \label{1}
 Q(F,\lambda) \ll
\lambda\int_{0}^{\lambda^{-1}} {|\widehat{F}(t)|
\,dt},\quad \lambda>0,
\end{equation}
where $\widehat{F}(t)$ is the corresponding characteristic
function.
 In the general case $Q(F,\lambda)$ cannot be estimated from below by the right hand side of
 inequality~\eqref{1}.
 However, if we assume additionally that the distribution $F$ is symmetric and
 its characterictic function is non-negative for all~$t\in\mathbf R$, then we have the lower bound:
\begin{equation} \label{1a}Q(F,\lambda)\gg
\lambda\int_{0}^{\lambda^{-1}} {\widehat{F}(t) \,dt}
\end{equation}
and, therefore,
\begin{equation} \label{1b}
Q(F,\lambda)\asymp\lambda\int_{0}^{\lambda^{-1}}
{\widehat{F}(t) \,dt}
\end{equation} (see \cite{Arak and Zaitsev}, Lemma~1.5 of Chapter II).
 The use of relation \eqref{1b} allows us to simplify the arguments of Friedland and Sodin
 \cite{Friedland and Sodin},  Rudelson and Vershynin~\cite{Rudelson and Vershynin} and Vershynin \cite{Vershynin} which were applied to the Littlewood--Offord problem
 (see also \cite{25, Eliseeva and Zaitsev}).
 \medskip

\emph{Proof of Theorem\/ $\ref{th1a}$.} Let $r$ be a fixed number satisfying $1<r\le\sqrt 2$.
Represent the distribution $G=\mathcal{L}(\widetilde{X})$
as a mixture
$$G=q E +\sum_{j=0}^{\infty}p_j G_j, $$ where $q={\mathbf
P}(\widetilde{X}=0)$, $p_j={\mathbf P}(\widetilde{X} \in A_j)$,
$j=0,1,2,\ldots$, $A_0=\{x:|x|>1\}$, $A_j=\{x: r^{-j}<|x|\leq
r^{-j+1}\}$,  $E$ is probability measure concentrated in zero,
$G_j$ are probability measures defined for $p_j>0$ by the formula
$$G_j\{X\}=G\{X\cap A_j\}/{p_j},$$ for any Borel
set~$X$.  In fact, $G_j$ is the conditional distribution of
$\widetilde X$ provided that $\widetilde X\in A_j$. If $p_j=0$,
then we can take as $G_j$ arbitrary measures.

For $z\in \mathbf{R}$, $\gamma>0$, introduce the distribution
$H_{z,\gamma}$, with the characteristic function
\begin{equation} \label{11}\widehat{H}_{z,\gamma}(t)=\exp\Big(-\cfrac{\gamma}{\,2\,}\;\sum_{k=1}^{n}\big(1-\cos(2a_k
zt)\big)\Big).\end{equation}
It is clear that $H_{z,\gamma}$ is a symmetric infinitely divisible distribution.
  Therefore, its characteristic function is positive for all $t\in \mathbf{R}$.

For the characteristic function $\widehat{F}(t)= \mathbf{E}\,\exp(itX)$, we have
$$|\widehat{F}(t)|^2 = \mathbf{E}\,\exp(it\widetilde{X}) =
\mathbf{E}\,\cos(t\widetilde{X}),$$
where $\widetilde{X}=X_1-X_2$ is the corresponding symmetrized random variable. Hence,
\begin{equation}\label{6}|\widehat{F}(t)| \leq
\exp\Big(-\cfrac{\,1\,}{2}\;\big(1-|\widehat{F}(t)|^2\big)\Big)  =
\exp\Big(-\cfrac{\,1\,}{2}\;\mathbf{E}\,\big(1-\cos(t\widetilde{X})\big)\Big).
\end{equation}

 According to \eqref{1} and \eqref{6}, we have
\begin{align}Q(F_a,1/D)&=Q(F_{2a},2/D)\le 2\,Q(F_{2a},1/D)\nonumber
\\ &\ll \frac 1D\int\limits_{0}^{D}|\widehat F_{2a}(t)|\,dt\nonumber
\\
 &\ll\frac 1D
\int\limits_{0}^{D}\exp\Big(-\frac{\,1\,}{2}\,\sum_{k=1}^{n}\mathbf{E}\,\big(1-\cos(2a_k
t \widetilde{X})\big)\Big)\,dt=I.\label{ww}
\end{align}
It is evident that
\begin{eqnarray*}
\sum_{k=1}^{n}\mathbf{E}\big(1-\cos(2a_k t
\widetilde{X})\big)&=&\sum_{k=1}^{n}\int_{-\infty}^{\infty}\big(1-\cos(2a_k
t x)\big)\,G\{dx\}
 \\
&=&\sum_{k=1}^{n}\sum_{j=0}^{\infty}\int_{-\infty}^{\infty}\big(1-\cos(2a_k
t x)\big)\,p_j
\,G_j\{dx\}\\
&=&\sum_{j=0}^{\infty}\sum_{k=1}^{n}\int_{-\infty}^{\infty}\big(1-\cos(2
a_k t x)\big)\,p_j \,G_j\{dx\}.
\end{eqnarray*}

 We denote $\beta_j=r^{-2j}p_j $,
$\beta=\sum_{j=0}^{\infty}\beta_j$, $\mu_j={\beta_j}/{\beta}$,
$j=0,1,2,\ldots$. It is clear that $\sum_{j=0}^{\infty}\mu_j=1$ and
${p_j}/{\mu_j}=r^{2j}\beta $ (for $p_j> 0$).

Let us estimate the quantity $\beta$:
\begin{eqnarray*}
\beta = \sum_{j=0}^{\infty}\beta_j &=&\sum_{j=0}^{\infty}r^{-2j} p_j  \,
= {\mathbf P}\big\{|\widetilde{X}|>1\big\} +
\sum_{j=1}^{\infty}r^{-2j}\,{\mathbf
P}\big\{r^{-j}<|\widetilde{X}|\leq r^{-j+1}\big\}  \\
&\geq&\int\limits_{|x|>1}\,G\{dx\} + \sum_{j=1}^{\infty}
\int\limits_{r^{-j}<|x|\leq r^{-j+1}}\cfrac{x^2}{r^{2}}\,G\{dx\}\\ &\geq&
\cfrac{\,1\,}{r^{2}}
\int\limits_{|x|>1}\,G\{dx\} + \cfrac{\,1\,}{r^{2}} \int\limits_{|x|\leq1}x^2 \,G\{dx\} =
\cfrac{\,1\,}{r^{2}} \,M(1).
\end{eqnarray*}
Since $1<r\le\sqrt 2$, this implies \begin{equation}\label{9}
\beta \geq \cfrac{\,1\,}{{2}}\, M(1).\end{equation} Condition
$L^2\geq 1/{M(1)}$ and inequality \eqref{9} give the bound
\begin{equation}\label{99}L^2\beta \geq \cfrac{\,1\,}{{2}}\,.
\end{equation}

 We now proceed similarly to the proof of
a result of Ess\'een \cite{Esseen} (see \cite{Petrov}, Lemma 4 of Chapter II).
Using the H\"older inequality, it is easy to see that
 \begin{equation}
 \label{66}
 I\leq \prod _{j=0}^{\infty}I_j^{\mu_j},
 \end{equation}
 where
\begin{eqnarray*}
I_j&=&\frac 1D\int_{0}^{D}\exp\Big(-\cfrac{p_j}{2\,\mu_j}\;\sum_{k=1}^{n}\int_{-\infty}^{\infty}\big(1-\cos(2a_k
t x)\big)\,G_j\{dx\}\Big)\,dt \\
&=&
\frac 1D\int_{0}^{D}\exp\Big(-\frac{\,1\,}2\,r^{2j}\beta\;\sum_{k=1}^{n}\int_{A_j}\big(1-\cos(2a_k
t
x)\big)\,G_j\{dx\}\Big)\,dt
\end{eqnarray*}
if $p_j > 0$, and
$I_j=1$ if $p_j=0$.

Applying Jensen's inequality to the exponential in the integral (see
\cite{Petrov},
p. 49)), we obtain
\begin{eqnarray}
I_j&\leq&\frac 1D\int_{0}^{D}\int_{A_j}\exp\Big(-\frac{\,1\,}2\,r^{2j}\beta\;
\sum_{k=1}^{n}\big(1-\cos(2a_k t x)\big)\Big)\,G_j\{dx\}\,dt\nonumber \\
&=
&\frac 1D\int_{A_j}\int_{0}^{D}\exp\Big(-\frac{\,1\,}2\,r^{2j}\beta\;\sum_{k=1}^{n}\big(1-\cos(2a_k
t x)\big)\Big)\,dt\,G_j\{dx\}\nonumber \\
&\leq& \sup_{z\in A_j}\frac 1D\int_{0}^{D}\widehat{H}_{z,1}^{r^{2j}\beta}(t)\,dt.\label{rr}
\end{eqnarray}

Let us estimate the characterictic function
$\widehat{H}_{\pi,1}(t)$ for $|\,t|\leq D$. We can proceed in the
same way as the authors of~\cite{Friedland and Sodin},
\cite{Rudelson and Vershynin} and \cite{Vershynin}. It is evident
that  $1-\cos x \geq 2x^2/\pi^2$, for~${|x|\leq\pi}$. For
arbitrary~$x$, this implies that
$$1-\cos x\geq 2\,\pi^{-2}
\min_{m\in \mathbf{Z}}|\,x-2\pi m|^2.$$ Substituting this
inequality into \eqref{11}, we obtain
\begin{eqnarray}
\widehat{H}_{\pi,1}(t)&\leq&\exp \Big(-\cfrac{1}{\pi^{2}} \;\sum_{k=1}^{n}\min_{m_k \in
\mathbf{Z}}\big|2\pi t a_k -2 \pi m_k\big|^2\Big)\nonumber \\
&=&\exp\Big(- 4\;\sum_{k=1}^{n}\min_{m_k \in
\mathbf{Z}}|\,ta_k-m_k|^2\Big)\nonumber \\
&=&\exp\big(- 4\;\big(\hbox{dist}(ta,\mathbf{Z}^n)\big)^2\big).\label{7b}
\end{eqnarray}
Using \eqref{4s}, wee see that, for $|\,t|\leq 1/2\,\|a\|_{\infty}$,
inequality \eqref{7b} turns into
\begin{equation}\label{7a}
\widehat{H}_{\pi,1}(t)\leq\exp(-4\,t^2).
\end{equation}

Now we can use relations \eqref{5bh},  \eqref{7b}
and \eqref{7a}  to estimate the integrals~$I_j$.
First we consider the case $j=1,2,\ldots$. Note that
the characteristic functions~$\widehat{H}_{z,\gamma}(t)$ satisfy the equalities
\begin{equation} \label{5}
\widehat{H}_{z,\gamma}(t)=\widehat{H}_{y,\gamma}\big({zt}/{y}\big)\quad\hbox{and}\quad
 \widehat{H}_{z,\gamma}(t)=\widehat{H}_{z,1}^{\gamma}(t).
\end{equation}
The first equality \eqref{5} implies that
\begin{equation} \label{55}
\hbox{if}\quad{H}_{z,\gamma}=\mathcal L(\xi),\quad\hbox{then}\quad {H}_{y,\gamma}=\mathcal L(y\,\xi/z).
\end{equation}

For $z\in A_j$ we have $r^{-j}<|z|\leq r^{-j+1}<\pi$. Hence, for
${|\,t|\leq D}$, we have $|{zt}/{\pi}|<D$. Therefore, using
properties \eqref{5} with $y=\pi$ and aforementioned estimates
\eqref{5bh},  \eqref{7b} and~\eqref{7a}, we obtain, for $z\in A_j$
and for $z=\pi$,
\begin{eqnarray*} \label{5cq}
\widehat{H}_{z,1}(t)
&\leq&\exp\big(-4\,f_L^2({zt}/{\pi})\big)\\
&=& \begin{cases}\qquad \exp\big(-({zt}/{\pi})^2/9\big),
& \hbox{ for }0<t\le eL\pi/z,\\ \exp\big(-4\,L^2\,\log(zt/L\pi)\big),&\hbox{ for }t> eL\pi/z.
\end{cases}
\end{eqnarray*}
Hence,
\begin{equation}\label{jj}
\sup_{z\in
A_j}\int\limits_{0}^{D}\widehat{H}_{z,1}^{r^{2j}\beta}(t)\,dt\leq
\int\limits_{0}^{D}\exp\big(-t^2\beta/9\pi^2\big)\,dt   +
\int\limits_{r^{j-1} L\pi
e}^{\infty}\Big(\frac{r^{j}L\pi}t\Big)^{4\,r^{2j}\beta L^2}\,dt
\ll \cfrac{1}{\sqrt{\beta}} \,.
\end{equation}
In the last inequality we used inequality \eqref{99}.

Consider now the case $j=0$. Relation \eqref{55} yields, for
$z>0,\,\gamma>0$,
\begin{equation}\label{8c}
Q(H_{z,\gamma},1/D)=Q\big(H_{1,\gamma},{1}/D{z}\big).
\end{equation}
Thus, according to  \eqref{8a}, \eqref{1b}, \eqref{5} and
\eqref{8c}, we obtain
\begin{eqnarray}
\sup_{z\in A_0}\frac 1D\int_{0}^{D}\widehat{H}_{z,1}^{\beta} (t) \,dt &=&
\sup_{z> 1} \frac 1D\int_{0}^{D}\widehat{H}_{z,\beta} (t) \,dt \asymp
\sup_{z> 1}\; Q(H_{z,\beta},1/D)\nonumber\\ &=&
\sup_{z> 1}\; Q\big(H_{1,\beta},{1}/{Dz}\big)
\leq Q(H_{1,\beta},1/D)\nonumber \\ &\asymp& Q\big(H_{1,\beta},{1}/D{\pi}\big) =
Q(H_{\pi,\beta},1/D)\nonumber\\ &\asymp&
\frac 1D\int_{0}^{D} \widehat{H}_{\pi,\beta}(t)\, dt =
\frac 1D\int_{0}^{D}\widehat{H}_{\pi,1}^{\beta}(t) \,dt.\label{pp}
\end{eqnarray}

Using the bounds \eqref{5bh}, \eqref{7b} and \eqref{7a} for the characteristic
function~$\widehat{H}_{\pi,1}(t)$ and taking into account inequality \eqref{99}, we have:
\begin{equation}\label{ee}
\int_{0}^{D}\widehat{H}_{\pi,1}^\beta(t)\,dt \leq
\int_{0}^{D}\exp(-t^2\beta/9)\,dt   +
\int_{Le}^{\infty}\Big(\frac{L}t\Big)^{4\beta L^2}\,dt \ll
\cfrac{1}{\sqrt{\beta}} \,.
\end{equation}

According to \eqref{rr}, \eqref{jj}, \eqref{pp} and \eqref{ee},
we obtained the same estimate
\begin{equation}\label{yy}
I_j \ll \cfrac{\,1\,}{D\sqrt{\beta}}\,,
\end{equation}
for all integrals $I_j$ with $p_j\neq 0$. In view of
$\sum_{j=0}^{\infty}\mu_j=1$, from \eqref{66} and \eqref{yy} it follows  that
\begin{equation}\label{yy8}I\leq\prod_{j=0}^{\infty}I_j^{\mu_j} \ll \cfrac{\,1\,}{D\sqrt{\beta}} \,.
\end{equation}
Using  \eqref{ww}, \eqref{9} and \eqref{yy8}, we complete the proof. $\square$

\medskip

Now we will deduce Corollary $\ref{c1b}$ from Theorem \ref{th1a}.
\medskip

\emph{Proof of Corollary $\ref{c1b}$.} We denote $b=a /\|a\|\in
\mathbf{R}^n$. Then the equality
$Q(F_a,\lambda)=Q(F_b,\lambda/\|a\|)$, for all $\lambda\ge0$,
holds. The vector~$b$ satisfies the conditions of Theorem
\ref{th1a} (which hold for the vector~$a$) with replacing $D$ by
$D\|a\|$. Indeed, $\|ub-m\|\geq f_L(u) $ for $u \in
\Big[\cfrac{1}{2\,\|b\|_{\infty}},D\|a\|\Big]$ and for all $m\in
\mathbf Z^n$. This follows from condition~\eqref{5bh} of Theorem
\ref{th1a}, if we denote $u={t}\|a\|$. It remains to apply
Theorem~\ref{th1a} to the vector $b$. $\square$
\bigskip

{\bf Acknowledgements.} The first and the third authors are supported by grant RFBR 10-01-00242.
The second and the third authors are supported by the SFB 701 in Bielefeld.
The third author is supported by grant RFBR 11-01-12104  and by the Program of
Fundamental Researches of Russian Academy of Sciences ``Modern
Problems of Fundamental Mathematics''.

\end{document}